\documentclass{article}
\title{Polynomial
and Rational Solutions of Holonomic Systems}
\author{Toshinori Oaku, Nobuki Takayama and Harrison Tsai}
\date{January 9, 1999}

\def\pd#1{ \partial_{#1} }
\def\qed{ {\tt [\kern-.4mm]} \bigbreak }
\def\k{ {\bf k} }
\def\span{\mbox{\rm Span}}
\newtheorem{theorem}{Theorem}[section]	
\newtheorem{proposition}[theorem]{Proposition}

\newtheorem{corollary}[theorem]{Corollary}

\newtheorem{remark}[theorem]{Remark}
\newtheorem{example}[theorem]{Example}
\newtheorem{algorithm}[theorem]{Algorithm}

\def\hom{\mbox{\rm Hom}}
\def\der{\mbox{\rm Der}}
\def\ext{\mbox{\rm Ext}}
\def\tor{\mbox{\rm Tor}}

\def\im{\mbox{\rm Im}}

\begin{document}

\maketitle

\section{Introduction}

Polynomial and rational solutions
for linear ordinary differential equations can
be obtained by algorithmic methods.
For instance, the maple package
{\tt DEtools} provides
efficient functions {\tt polysols} and {\tt ratsols}
to find polynomial and rational solutions for a given
linear ordinary differential equation with rational function
coefficients.

A natural analogue of the notion of linear ordinary 
differential equation in the several
variable case is the notion of holonomic system.
A holonomic system is a system of linear partial differential equations
whose characteristic variety is middle dimensional.

Chyzak~\cite{Chyzak} gave an algorithm to find the rational
solutions of holonomic systems by using elimination in the ring
of differential operators with rational function coefficients
combined with Abramov's algorithm for rational solutions of
ordinary differential equations with parameters.
To the authors,
solving holonomic systems is analogous to solving systems
of algebraic equations of zero-dimensional ideals.
Under this analogy, the method of Chyzak corresponds 
to the elimination method for solving systems of algebraic 
equations.

The aim of this paper is to give two new algorithms,
which are elimination free, 
to find polynomial and rational solutions 
for a given holonomic system associated to
a set of linear differential operators in the Weyl algebra
$$D = \k \langle x_1, \ldots, x_n, \partial_1, \ldots, \partial_n \rangle$$
where $\k$ is  a subfield of ${\bf C}$.

Polynomial and rational solutions can be obtained,
if they exist, 
by using an exhaustive search.
For instance, when $f=0$ is the singular locus of a holonomic system
$M = D/I$,
any rational solution has the form
$ g/f^r$.
If we have upper bounds for the degree of the polynomial $g$ and for $r$,
then we can construct all rational solutions
by solving linear equations satisfied by the coefficients of $g$.
Alternatively, if we know the dimension of rational solutions,
then we can obtain all rational solutions by increasing the degree
of $g$ and $r$.
Hence, the problem reduces to finding effective bounds for these
numbers.

In sections 2 and 3, we give algorithms for upper bounds on
the degree of $g$ and on $r$.  The main techniques we use
are Gr\"obner deformations in $D$ as introduced in the book
\cite{SST-BOOK} and the $b$-function
for $D/I$ and $f$.

In section 4, we give an algorithm to evaluate the dimension
of polynomial and rational solutions.
Our approach is an analog in $D$ of a question studied
by Singer~\cite{Singer}, who
gave an algorithm to compute ${\rm Hom}_R(M,N)$
for left $R := \k(x_1)\langle \partial_1 \rangle$-modules
$M$ and $N$ and studied its relation to factorizations of ordinary
differential operators.
The theory of $D$-modules translates our problem 
on polynomial and rational solutions
to constructions in the ring of differential operators $D$.
For example, the $\k$-vector space
$$
  {\rm Hom}_D (D/I, \k[{\bf x}]) \simeq H^{-n}(\Omega \otimes_D^L {\bf D}(D/I))
$$
is the space of the polynomial solutions of the left ideal $I$.
Here, $\Omega$ is the module of the top dimensional differential forms
and ${\bf D}$ is the dualizing functor.
See, e.g., the book of Bj\"ork \cite{Bj} on this translation.
We evaluate the dimension of the right hand side by
recent developments of computational algebra such as 
construction of free resolutions in the ring $D$ and restrictions
of $D$-modules \cite{O}, \cite{OT}, \cite{Ts}, \cite{W1}.
Our method also allows us to evaluate the dimension of solutions to
a holonomic system inside any holonomic module.  For instance, we can 
find the dimension of the delta function solutions to $I$.

Throughout the paper, we refer to the book \cite{SST-BOOK}
for fundamental facts on the algorithmic treatment of $D$.
Also, the algorithms which appear in the paper have
been implemented in either {\tt kan}~\cite{kan} or 
Macaulay~2~\cite{Macaulay2}.

We deeply thank Dan Grayson, Anton Leykin, and Mike Stillman,
who helped to implement $D$-modules in Macaulay~2, and
Fr\'ed\'eric Chyzak and Michael Singer for discussions
on rational solutions.

\section{Polynomial solutions by Gr\"obner deformations}
\label{GBdeform}

How can we obtain all polynomial solutions for ordinary differential
equations?
One method is to compute the indicial polynomial at infinity,
find an upper bound on the degrees of polynomial solutions,
and determine the coefficients of polynomials.
The analogous method works for holonomic systems by using Gr\"obner
deformations.
For $\ell \in D$ and the weight vector $w \in {\bf R}^n$,
we denote by 
${\rm in}_{(-w,w)}(\ell)$
the initial term of $\ell$ with respect to the weight
$(-w,w)$
(see, e.g., \cite[\S 1.1]{SST-BOOK}).  The following
proposition follows from the definition of ${\rm in}_{(-w,w)}(\ell)$.

\begin{proposition}
\label{vertex}
Suppose that $f(x_1, \ldots, x_n)$ is a polynomial solution
of 
$I = D \cdot \{ \ell_1, \ldots, \ell_m \}$.
Take $w \in {\bf Z}^n$.
Then
$ f(t^{w_1}x_1, \ldots, t^{w_n} x_n)$ can be expanded as a polynomial in $t$
as
$$  f_w({\bf x}) t^p + O(t^{p+1}). $$
Then we have
$$  {\rm in}_{(-w,w)}(\ell_i) \bullet f_w = 0. $$
\end{proposition}

The initial ideal
${\rm in}_{(-w,w)}(I)$ is sometimes called
the Gr\"obner deformation of
$I$ with respect to $(-w,w)$. 

\begin{theorem} {\rm \cite{Assi-Castro-Granger}}
There exist only finitely many Gr\"obner deformations. 
\end{theorem}

The Newton polytope of a polynomial solution $f$ is defined as
the convex hull of the exponent vectors of $f$.
For generic $w$, $f_w$ is a monomial $c x^a$ and the point $a$
is a vertex of the Newton polytope of $f$.

Let $R = \k ( x_1, \ldots, x_n) \langle \pd{1}, \ldots, \pd{n} \rangle$
and $\theta_i = x_i \pd{i}$.
Since $a \in {\bf Z}^n$ belongs to the zero set of the indicial ideal
$$ \widetilde{{\rm in}_{(-w,w)}}(I) =
  R \cdot {\rm in}_{(-w,w)}(I) \cap \k [ \theta_1, \ldots, \theta_n] ,
$$
we can construct a polytope that contains the Newton polytopes
of the polynomial solutions by taking the convex hull of
all the non-negative integral roots of all possible indicial ideals.

It is not necessary to find all Gr\"obner deformations
to obtain polynomial solutions.
Let $b(s)$ be the generator of
${\rm in}_{(-w,w)}(I) \cap \k [ s ] $,
$ s = \sum_{i=1}^n w_i \theta_i$.
The polynomial $b(s)$ is called the $b$-function of $I$
with respect to $(-w,w)$.
The next proposition follows from the definition of $b(s)$.
\begin{proposition}
Let $w$ be a strictly negative weight vector. In other words,
we assume that $w_i < 0$ for all $i$.
Consider the b-function $b(s)$ of $I$ with respect to $(-w,w)$
and let $-k_1$ be the smallest integer root of $b(s)=0$.
The polynomial solutions of $I$ have the form
\begin{equation} \label{eq:shape}
 \sum_{ p_i \geq 0,\, p\cdot w \leq k_1 } c_p x^p.
\end{equation}
\end{proposition}

\begin{algorithm}
\label{alg-polygrob}
\rm
(Finding the polynomial solutions by a Gr\"obner deformation)

\noindent{\sc Input}: a holonomic left ideal $I$.
 
\noindent{\sc Output}: the polynomial solutions of $I$.

\begin{enumerate}
\item Take a strictly negative weight vector $w$,
compute the Gr\"obner deformation ${\rm in}_{(-w,w)}(I)$,
and compute the smallest non-positive
integer root $-k_1$ of the $b$-function 
with respect to $(-w,w)$.
See, e.g., \cite[Alg. 5.15]{SST-BOOK} for these procedures.
\item If we do not have such a root, then there is no polynomial solution
other than $0$.
\item If there is a minimal integer root,
then determine the coefficients $c_p$ of (\ref{eq:shape})
by solving linear equations for the coefficients.
\end{enumerate}
\end{algorithm}

\begin{example}
\label{appell} \rm
The following system of differential equations of two variables
is called 
the Appell differential equation $F_1(a,b,b',c)$ \cite{Appell}:
\begin{eqnarray*}
&&\theta_x ( \theta_x+\theta_y+c-1) - x ( \theta_x+\theta_y+a) (\theta_x+b), \\
&&\theta_y ( \theta_x+\theta_y+c-1) - y ( \theta_x+\theta_y+a) (\theta_y+b'), \\
&& (x-y) \pd{x}\pd{y} - b' \pd{x} + b \pd{y} \\
\end{eqnarray*}
where $a, b, b', c$ are complex parameters.
Let us demonstrate how Algorithm \ref{alg-polygrob} works
for the system of parameter values
$ (a,b,b',c)= (2, -3, -2, 5)$.
First, we choose a strictly negative weight vector $w =(-1,-2)$ and
compute the $b$-function $b(s)$, $s = -\theta_x - 2 \theta_y$, 
which is the generator of the principal ideal
$ {\rm in}_{(-w,w)}(I) \cap {\bf Q}[-\theta_x - 2 \theta_y]$.
We can use the V-homogenization or the homogenized Weyl algebra
to get the generator
(see, e.g., \cite[\S 1.2]{SST-BOOK}).
Second, we need to find the integer roots of the $b$-function $b(s)=0$.
In our example,
these are 
$$ -7, 0, 4. $$
From Proposition \ref{vertex}, 
the highest $(-w)$-degree monomial $ c x^p y^q$
in a polynomial solution gives rise to an integer solution 
$ w_1 p + w_2 q = -p - 2q $ of the $b$-function.
Hence, the polynomial solutions are of the form
$$  f = \sum_{p,q \geq 0, p + 2q \leq 7} c_{pq} x^p y^q. $$
Finally, 
we determine the coefficients $c_{pq}$
by applying the differential operators to $f$ and putting the results
to $0$.
In our example, 
we have only one polynomial solution
\begin{eqnarray*}
&& (-\frac{1}{21}y^2+\frac{1}{7} y-\frac{4}{35}) x^3+
   (\frac{3}{14} y^2-\frac{24}{35} y+\frac{3}{5}) x^2 \\
&+&(-\frac{12}{35} y^2+\frac{6}{5} y-\frac{6}{5}) x
    +\frac{1}{5} y^2-\frac{4}{5} y+1.
\end{eqnarray*}
\end{example}

\section{Rational solutions by Gr\"obner deformations}

The singular locus of a $D$-ideal $I$ is defined to be
the projection of the characteristic variety of $I$ minus the zero section
from the cotangent bundle to the coordinate base space.  In other words,
it is the zero set
$${\rm Sing}(I) = V\left( \langle 
{\rm in}_{(0,e)}(I) : (\xi_1,\dots,\xi_n)^{\infty}\rangle \cap
\k[x_1,\dots,x_n] \right).$$

Any rational solution to $I$ has its poles contained inside
the singular locus.
Thus if $f({\bf x})$ defines the codimension 1 component
of ${\rm Sing}(I)$, we may limit our search for
rational solutions to $\k[{\bf x}][\frac{1}{f}]$.

We will present a method to obtain an upper bound of the order 
of the poles along $f=0$ for each rational solution. 
For this purpose we use the notion of the $b$-function for
$f$ and a section $u$ of a holonomic system, which was introduced 
by Kashiwara \cite{KashiwaraHolonomicII}:
Let ${\cal D}$ be the sheaf of algebraic differential operators
on $X = {\bf C}^n$.
For a holonomic ${\cal D}$-module ${\cal M}={\cal D}/{\cal D}I$ 
and a polynomial $f$,
consider the tensor product
\begin{equation}
\label{eqn-bfunction}
{\cal N} = {\cal O}[f^{-1},s]f^s \otimes_{{\cal O}_X} 
{\cal M}.
\end{equation}
This ${\cal N}$ has a structure of a left ${\cal D}$-module via 
the Leibnitz rule.
Let $u$ be a section of ${\cal M}$. 
Then the  $b$-function for $f$ and $u$ (or for $f^su$) 
at $p \in {\bf C}^n$ is the minimum degree monic polynomial 
$0 \not= b(s) \in {\bf C}[s]$ such that
\begin{equation}
\label{eqn-bfunction2}
b(s)f^s\otimes u \in {\cal D}[s](f^{s+1} \otimes u)
\end{equation}
holds in ${\cal N}$ at $p$ (i.e., as a germ of ${\cal N}$ at $p$). 
This $b$-function depends on the point $p$.  As a function of $p$,
there is a stratification of ${\bf C}^n$ for which
the $b$-function does not change on each strata (see e.g.~\cite{O}
for an algorithmic proof of this fact).  In the
definitions~(\ref{eqn-bfunction}) and~(\ref{eqn-bfunction2}) for
$b$-function,
if we replace ${\cal O}$ by the polynomial ring ${\bf k}[{\bf x}]$, 
${\cal D}$ by the Weyl algebra $D$,
and ${\cal M}$ by a holonomic $D$-module $M = D/I$,
then we obtain the global $b$-function for $f$ and $u$.
It is the least common multiple of $b$-functions at every point.

\begin{theorem}  \label{proposition:upperbound}
Let $u$ be the residue class of $1$ in ${\cal D}/{\cal D}I$,
and let $b(s)$ be the $b$-function for $f$ and $u$
at a point $p \in {\bf C}^n$ where $f(p) = 0$.
Assume that $I$ admits an analytic solution of the form $gf^r$ 
around $p$, where $r \in {\bf C}$, $g$
is a holomorphic function on a neighborhood of $p$, 
and $g(p) \neq 0$. Then $s+r+1$ divides $b(s)$.
\end{theorem}

\noindent \textbf{Proof}:
Let ${\cal D}^{\rm an}$ and 
${\cal O}^{\rm an}$ be respectively the sheaf of analytic differential operators
and the sheaf of holomorphic functions on ${\bf C}^n$.
We may define the analytic $b$-function by replacing
${\cal O}$ by ${\cal O}^{\rm an}$,
${\cal D}$ by ${\cal D}^{\rm an}$, and
${\cal M}$ by a ${\cal D}^{\rm an}$-module ${\cal M}^{\rm an}$
in the definitions~(\ref{eqn-bfunction}) 
and~(\ref{eqn-bfunction2}).
Since the $b$-function is an analytic invariant and the analytic and the 
algebraic $b$-functions
coincide (see e.g. \cite[\S 8]{O}), we may work in 
the analytic category.  We do this to consider solutions $gf^r$ where $g$
is holomorphic at $p$.  If we only wish to consider solutions $gf^r$ where 
$g$ is a polynomial, then we may work in the algebraic category.

In general, given a map of left ${\cal D}^{\rm an}$-modules 
$\phi:{\cal M}_1^{\rm an}
\rightarrow {\cal M}_2^{\rm an}$ and a section $u$ of ${\cal M}_1^{\rm an}$, 
the $b$-function for $f^su$ at a point $p$
is divisible by the $b$-function for $f^s\phi(u)$ at $p$.
We apply this basic fact to the following map $\varphi$.
Let $J^{\rm an}$ be the annihilating ideal of $g f^r$ in 
${\cal D}^{\rm an}$.
Since $J^{\rm an} \supseteq I^{\rm an}
:= {\cal D}^{\rm an}I$ and $g(p) \neq 0$, 
we have a left ${\cal D}^{\rm an}$-homomorphism
$$
\varphi : {\cal D}^{\rm an}/I^{\rm an} 
\longrightarrow {\cal D}^{\rm an}gf^r = {\cal D}^{\rm an}f^r
\hookrightarrow {\cal O}^{\rm an}[f^{-1}]f^r
$$
which sends $u$ to $g f^r$.  
This map extends to a left ${\cal D}^{\rm an}[s]$-homomorphism
$$\begin{array}{cl}
1\otimes \varphi : &
{\cal O}^{\rm an}[f^{-1},s]f^s
\otimes_{{\cal O}^{\rm an}}
{\cal D}^{\rm an}/I^{\rm an} \\
 & \\
& \longrightarrow
{\cal O}^{\rm an}[f^{-1},s]f^s
\otimes_{{\cal O}^{\rm an}}
{\cal O}^{\rm an}[f^{-1}]f^r
= {\cal O}^{\rm an}[f^{-1},s]f^{s+r}
\end{array}$$
which sends $f^s \otimes u$ to $gf^{s+r}$.
By the definition of $b(s)$, there exists a germ $P(s)$ of
${\cal D}[s]$ at $p$ such that 
$$ (P(s)f - b(s))(f^s \otimes u) = 0. $$
Since $1\otimes \varphi$ is a left ${\cal D}^{\rm an}$-homomorphism,
applying it to the above equation gives the equation
$(P(s)f - b(s))(gf^{s+r}) = 0$,
or in other words, 
$$g^{-1}P(s)gf^{s+r+1} = b(s)f^{s+r}.$$
Thus, we see that the Bernstein-Sato polynomial $b_f(s)$ of $f$ at $p$
divides $b(s-r)$.  
Note that $s+1$ divides $b_f(s)$ since $f(p) = 0$ 
(cf.\ \cite{KashiwaraBfunction}).
In conclusion, we have proved that $s+1$ divides $b(s-r)$. 
This completes the proof.  
\qed

By virtue of the above theorem, we can obtain upper bounds
by computing the $b$-function for $f^su$ at a smooth point of each 
irreducible component of the singular locus of $I$.  From now on,
let us also take $f \in \k[{\bf x}]$ to be a square-free
polynomial defining the codimension one component of the singular
locus, and let $f = f_1\cdots f_m$ be its irreducible decomposition
in $\k[{\bf x}]$.

\begin{theorem}  \label{theorem:upperbound}
Let $b_i(s)$ be the $b$-function for $f_i^su$ at a generic point
of $f_i = 0$.  
Denote by $r_i$ the maximum integer root of $b_i(s)=0$.
Then any rational solution (if any) to $I$ can be written in the form
$g f_1^{-r_1-1}\cdots f_m^{-r_m-1}$ with a polynomial 
$g \in {\bf C}[{\bf x}]$.
If some $b_i(s)$ has no integral root, then there exist no 
rational solutions to $I$ other than zero.
\end{theorem}

\noindent \textbf{Proof:}
An arbitrary rational solution to $M$ is written in the form
$gf_1^{-\nu_1}\cdots f_m^{-\nu_m}$ with integers $\nu_1,\dots,\nu_m$ 
and $g \in {\bf C}[{\bf x}]$.
Since the space of the rational solutions with coefficients in 
${\bf C}$ is spanned by those with coefficients in $\k$, 
we may assume $g \in \k[x]$, and $f$ and $g$ are relatively prime
in $\k[{\bf x}]$. 
Let $p$ be a generic point of $f_i = 0$. 
We may assume that $f_i$ is smooth at $p$, $g(p) \neq 0$,
and $f_j(p) \neq 0$ for $j \neq i$.
It follows from Theorem \ref{proposition:upperbound} 
that $b_i(\nu_i-1) = 0$.  This implies $\nu_i \leq r_i+1$.
\qed

Since $b$-functions divide the global $b$-function,
an upper bound can also be obtained from the global $b$-function.

\begin{corollary}
\label{cor:ratbound}
Let $b_i(s)$ be the global $b$-function for $f_i^su$,  
and denote by $r_i$ the maximum integer root of $b_i(s)=0$.
Then any rational solution (if any) to $I$ can be written in the form
$g f_1^{-r_1-1}\cdots f_m^{-r_m-1}$ with a polynomial 
$g \in {\bf C}[{\bf x}]$.
\end{corollary}

We mention the corollary 
since the algorithm to compute global
$b$-functions is simpler than the algorithm
to compute $b$-functions.  However, the $b$-function
offers finer information.  For instance, the well-known example
$f = x^2+y^2+z^2+w^2$ has Bernstein-Sato polynomial $(s+1)(s+2)$ coming
from the functional equation,
$\frac{1}{4}(\partial_x^2+\partial_y^2+\partial_z^2+\partial_w^2)
\cdot f^{s+1} = (s+1)(s+2)f^s$.
Now consider the module $M = D\cdot f^{-1}$ and
let $u$ be the section of $f^{-1}$.
The global $b$-function for $f^su$ is $s(s+1)$, and hence
Corollary~\ref{cor:ratbound} implies that
rational solutions of $M$ all have the form $gf^{-1}$
or $gf^{0}$, where $g$ is a polynomial not divisible by $f$.  
On the other hand, the Bernstein-Sato polynomial
of $f$ at any nonsingular point $p$ of $f = 0$ (i.e.~except
for the origin) is $s+1$. It follows that
the $b$-function for $f^su$ equals $s$ at the generic
point of $f=0$ and hence Theorem~\ref{theorem:upperbound}
implies that all rational solutions actually have
the form $gf^{-1}$.

An algorithm to compute the $b$-function and the global 
$b$-function for $f^su$ was first given in \cite{O} based upon 
tensor product computation, which is slow and memory intensive. 
Shortly thereafter, Walther introduced in \cite{W0}  
a more efficient method to compute the global $b$-function for $f^su$.  
Both methods give the global $b$-function exactly,
under the condition that $I$ is $f$-saturated.
Otherwise, we get a multiple of the global $b$-function.  
Similarly, the method of~\cite{O} gives the $b$-function
exactly if $I$ is $f$-saturated and additionally
a certain primary decomposition in ${\bf C}[{\bf x}]$
is known.  If primary decomposition is
only available in ${\bf k}[{\bf x}]$, we again
get a multiple of the $b$-function.

Let us now describe an algorithm to compute the $b$-function for 
$f^s u$ at a generic point of $f = 0$ by combining the method of \cite{W0} 
and the primary decomposition as was used in \cite{O}.  

\begin{algorithm} \label{algorithm:upperbound} \rm
(Computing an upper bound of the $b$-function at a generic point)

\noindent{\sc Input}: a finite set $G_0$ of generators of 
a holonomic $D$-ideal $I$ 
and an irreducible polynomial $f \in \k[{\bf x}]$.

\noindent{\sc Output}: $b'(s) \in \k[s]$, which is a multiple 
of the $b$-function $b(s)$ for $f^su$ at a generic point 
of $f = 0$, where $u$ is the residue class of $1$ in $D/I$. 

\begin{enumerate}
\item 
Introducing a new variable $t$, 
put $\vartheta_i = \pd{i} + (\partial f/\partial x_i)\pd{t}$.
Let $\widetilde I$ be the left ideal of $D_{n+1}$, the Weyl algebra on 
the variables $x_1,\dots,x_n,t$, that is generated by 
$$
\{ P({\bf x},\vartheta_1,\dots,\vartheta_n) \mid 
P({\bf x},\pd{1},\dots,\pd{n}) \in G_0 \} \cup \{ t - f({\bf x}) \}.
$$
\item
Let $G_1$ be a finite set of generators of the left ideal  
${\rm in}_{(-1,0,\dots,0;1,0,\dots,0)}(\widetilde I)$ of $D_{n+1}$.
Here, $-1$ is the weight for $t$ and $1$ is the weight for $\pd{t}$.
\item
Rewrite each element $P$ of $G_1$ in the form
$$
P = \pd{t}^\mu P'(t\pd{t},{\bf x},\pd{1},\dots,\pd{n})
\quad\mbox{or}\quad
P = t^\mu P'(t\pd{t},{\bf x},\pd{1},\dots,\pd{n})
$$
with a non-negative integer $\mu$, and define $\psi(P)$ by,
$$
\psi(P) := t^{\mu}\pd{t}^{\mu}P' = 
t\pd{t}\cdots(t\pd{t}-\mu+1)P'(t\pd{t},{\bf x},\pd{1},\dots,\pd{n})
$$
$$
\quad\mbox{or}\quad
\psi(P) := \pd{t}^{\mu}t^{\mu}P' = 
(t\pd{t}+1)\cdots(t\pd{t}+\mu)P'(t\pd{t},{\bf x},\pd{1},\dots,\pd{n}).
$$
Put
$$
G_2 := \{ \psi(P)(-s-1,{\bf x},\pd{1},\dots,\pd{n}) \mid P \in G_1\}.
$$
\item
Compute the elimination ideal $J := \k[s,{\bf x}] \cap D[s]G_2$.
(The global $b$-function can be obtained at this stage
by computing the monic generator of the ideal $J \cap \k[s]$.)
\item
Compute a primary decomposition of $J$ in $\k[s,{\bf x}]$ as
$$
J = Q_1 \cap \cdots \cap Q_\nu.
$$
\item
For each $i=1,\dots,\nu$, 
compute $Q_{ix} := Q_i \cap \k[{\bf x}]$, which is a primary ideal of $\k[{\bf x}]$.
\item
Let $b'(s)$ be the monic generator of the ideal 
$$
\bigcap \{ Q_i \cap \k[s] \mid \sqrt{Q_{ix}} \subset \k[{\bf x}]f \}
$$
of $\k[s]$. (Note that $\sqrt{Q_{ix}} \subset \k[{\bf x}]f$ implies 
that $\sqrt{Q_{ix}}$ equals $\k[{\bf x}]f$ or $\{0\}$.)
\end{enumerate}
\end{algorithm}

\begin{theorem} \label{proposition:indicial}
In the above algorithm, the polynomial
$b'(s)$ is precisely the $b$-function 
for $f^su$ at a generic point of $f=0$ 
if $I$ is $f$-saturated (i.e., $I:f^\infty = I$)
and each ${\bf C}[s,{\bf x}]Q_i$ remains primary in ${\bf C}[s,{\bf x}]$. 
Otherwise, the polynomial $b'(s)$ is a multiple of the $b$-function 
for $f^su$ at a generic point of $f=0$.
\end{theorem}

\noindent \textbf{Proof:}
Using essentially the same method as the proof of Lemma 4.1 in \cite{W0}, 
we can prove that $\widetilde I$ is precisely the annihilator ideal for 
$\delta(t-f({\bf x}))\otimes u$ in 
$$
\widetilde M := (D_{n+1}\delta(t-f({\bf x})))\otimes_{{\bf C}[{\bf x}]}D/I,
$$
where $\delta(t-f({\bf x}))$ denotes the residue class of $(t-f({\bf x}))^{-1}$ in
$\k[{\bf x},(t-f({\bf x}))^{-1}]/\k[{\bf x}]$.
Let $b_t(s)$ be the indicial polynomial for 
$\delta(t-f({\bf x}))\otimes u$ along $t=0$ at a point $(0,p)$ with 
$f(p) = 0$.
Then by Theorem 6.14 of \cite{O}, the $b$-function $b(s)$ for
$f^su$ at $p$ divides out, and if $I$ is $f$-saturated, coincides with 
$b_t(-s-1)$.

It follows from the definition that $b_t(-s-1)$ is a generator of the ideal
${\cal O}_p[s]J \cap {\bf C}[s]$ of ${\bf C}[s]$, where ${\cal O}_p$ 
denotes the stalk of ${\cal O}$ at $p$. 
If ${\bf C}[s,{\bf x}]Q_i$ are primary in ${\bf C}[s,{\bf x}]$, 
$b'(s)$ generates the above ideal 
in view of Theorem 4.7 of \cite{O} (cf.\ also Lemma 4.4 of \cite{OT}).
In general, although $Q_i$ is primary in $\k[s,{\bf x}]$,
the extension ${\bf C}[s,{\bf x}]Q_i$ is no longer primary 
in ${\bf C}[s,{\bf x}]$ and admits a primary 
decomposition 
$$ 
{\bf C}[s,{\bf x}]Q_i = Q_{i1} \cap \cdots \cap Q_{i\mu_i}.
$$
In this case,  $b'(s)$ is the least common multiple of the generators of
the ideals ${\cal O}_p[s]Q_{ij} \cap {\bf C}[s]$ for $j=1,\dots,\mu_i$, 
while $b_t(-s-1)$ is the generator of ${\cal O}_p[s]Q_{ij} \cap {\bf C}[s]$ 
for some $j$ (such that $p$ belongs to the zero set of 
$Q_{ij} \cap {\bf C}[{\bf x}]$ which is also the zero set of a factor of 
$f$).  This completes the proof.
\qed

\begin{remark} \rm
In the notation in the above proof, the linear factors of $b'(s)$ and 
those of $b_t(-s-1)$ coincide.  In particular the set of integer roots 
of $b'(s) = 0$ is the same as that of $b_t(-s-1)=0$.
In fact, this follows from the fact that the linear factors of $b'(s)$ 
in $\k[s]$ are invariant under the action of the Galois group of 
$\overline \k$ over $\k$. 
\end{remark}

\begin{remark} \rm
$I$ is $f$-saturated if and only if the $-1$-th cohomology group
of the restriction of $\widetilde M$ in the proof of Proposition 
\ref{proposition:indicial} to $t=0$ vanishes
(see Theorem 6.4 and Proposition 6.13 of \cite{O}), 
which is computable by Algorithm 5.10 of \cite{O}
or by Algorithm 5.4 of \cite{OT}.
\end{remark}

\begin{remark} \rm
The $f$-saturation of $I$, which is the ideal 
$D[\frac{1}{f}]\cdot I \cap D$,
may be computed by using the localization algorithm of~\cite{OTW}
(if $I$ is specializable along $f$) or by using the 
less efficient algorithm of~\cite{Ts} (if $I$ is general).
By replacing $I$ with its $f$-saturation, we may then compute
the local $b$-function exactly.  However, since saturation is often an
expensive algorithm, we avoid making this replacement
in practice.
\end{remark}

Once we have determined the integers $r_1,\dots,r_m$ of 
Theorem \ref{theorem:upperbound}, we can use Gr\"obner deformations
to obtain the rational solutions.
Put $k_i = r_i+1$.  Then by virtue of Theorem \ref{theorem:upperbound},
we have only to determine rational solutions of the form
$g f_1^{-k_1}\cdots f_m^{-k_m}$ for some polynomial $g$. 
This amounts to computing polynomial solutions of some twisted ideal 
$I_{(k_1,\dots,k_m)}$ of $I$.
Namely, consider how $f \partial_i$ acts on an element 
$g f_1^{-k_1}\cdots f_m^{-k_m}$:
$$
f \partial_i \bullet (g f_1^{-k_1}\cdots f_m^{-k_m}) =
\left(f\frac{\partial g}{\partial x_i} - 
\sum_{j=1}^m k_j\frac{f}{f_j}\frac{\partial f_j}{\partial x_i}g\right)
f_1^{-k_1}\cdots f_m^{-k_m}.
$$
In other words, $f\partial_i$ acts on the numerator $g$ as the
differential operator 
\begin{equation}
\label{twistsub}
L_i = f \partial_i - \sum_{j=1}^m k_j\frac{f}{f_j}
\frac{\partial f_j}{\partial x_i}.
\end{equation}
Thus, if we multiply a set of generators $\{g_1,\dots, g_m\}$
of $I$ by sufficiently high powers $f^{m_i}$ such that
$f^{m_i}g_i \in 
\k\langle x_1,\dots,x_n,f\partial_1,\dots,f\partial_n\rangle,$
and if $I'$ is the ideal obtained from
$\{f^{m_1}g_1,\dots,f^{m_r}g_r\}$ by substituting $L_i$ for
$f\partial_i$, then the rational solutions of $I$
have the polynomial solutions of $I'$ as numerators.
Thus, it only remains to compute the polynomial solutions of $I'$.
However, since $I'$ might not be holonomic, we cannot
apply Algorithm~\ref{alg-polygrob} just yet.

Hence let us define $I_{(k_1,\dots,k_m)} :=
\k({\bf x})\langle \partial \rangle \cdot I' \cap D$, which
is the Weyl closure of $I'$ and whose polynomial
solutions are the same as $I$.  The advantage of
$I_{(k_1,\dots,k_m)}$ is that it is indeed holonomic, 
which follows from a theorem of Kashiwara.  
Namely, note that since $I$ is of finite rank, $I'$
remains of finite rank because an element in $I$ of the form
$(g_i({\bf x})\partial_i^{N_i} + \mbox{ lower order elements})$
will be sent to an element in $I'$ of the form
$(f^{M_i}g_i({\bf x}) \partial_i^{N_i} + \mbox{ lower order elements})$.
Now let $h({\bf x})$ be any polynomial vanishing on
the singular locus of $I'$.
Then the non-holonomic locus of $I'$ is
contained inside the zero set of $h({\bf x})$ 
regarded as a function on the cotangent bundle, 
and a theorem due to Kashiwara~\cite{KashiwaraHolonomicII} 
states that the ideal $D[h^{-1}]\cdot I' \cap D$ is holonomic.  
Furthermore, an argument in~\cite{Ts} shows that the Weyl closure of $I'$
also equals $D[h^{-1}]\cdot I' \cap D$. 
Summing up, we arrive at the following algorithm.

\begin{algorithm}\rm
(Computing the rational solutions of a holonomic ideal)

\noindent{\sc Input}: generators of a holonomic $D$-ideal $I$.

\noindent{\sc Output}: A basis of the rational solutions 
$h \in \k({\bf x})$ of $I \bullet h = 0$.

\begin{enumerate}
\item 
Compute a polynomial $f$ defining the 
codimension 1 component of ${\rm Sing}(I)$.
\item
Compute the irreducible decomposition $f = f_1\cdots f_m$ in $\k[{\bf x}]$.
\item 
For each $i=1,\dots,m$, 
compute the output $b'(s)$ of Algorithm \ref{algorithm:upperbound} 
with $I$ and $f_i$ as input.
Let $r_i$ be the maximum integer root of $b'(s) = 0$ and put $k_i = r_i+1$.
If $b'(s)$ has no integral root for some $i$, then there exists no 
rational solution other than zero. 
\item Compute the twisted ideal $I_{(k_1,\dots,k_m)}$ as follows.
First, form the ideal $I'$ described in the paragraphs
preceding the algorithm.
Second, compute any polynomial $h({\bf x})$ vanishing
on the singular locus of $I'$.
Third, compute the localization $(D/I')[h^{-1}]$ using
the algorithm in~\cite{OTW}.  Then the ideal $I_{(k_1,\dots,k_m)}$
is the kernel of $D \rightarrow D/I' \rightarrow (D/I')[h^{-1}]$.
\item Compute a basis $\{g_1,\dots,g_k\}$ of the polynomial solutions 
of $I_{(k_1,\dots,k_m)}$ using Algorithm~\ref{alg-polygrob}.
\item Output: $\{g_1f_1^{-k_1}\cdots f_m^{-k_m},\dots,
g_kf_1^{-k_1}\cdots f_m^{-k_m}\}$, a basis of the rational solutions 
of $I$.
\end{enumerate}
\end{algorithm}

\begin{example} \rm
Let $I$ be the left ideal generated by
$$ L_1 = \theta_x (\theta_x+\theta_y) -x (\theta_x+\theta_y+3)(\theta_x-1) $$
and
$$ L_2 = \theta_y (\theta_x+\theta_y) -y (\theta_x+\theta_y+3)(\theta_y+1). $$
The Appell function $F_1(3,-1,1,1;x,y)$ is a solution of this system.
The singular locus of $I$ is $xy(x-y)(1-x)(1-y) = 0$. 
We can compute the local indicial polynomial of $u$, 
the modulo class of $1$ in $D_2/I$, along $x=0$ directly by the 
algorithm of \cite[Section 4]{O}:
It is $s(s-1)$ on $\{(0,y) \mid y \neq 0\}$, and $s(s-1)^2$ 
at $(0,0)$.  
In the same way, the indicial polynomial of $u$ along $y=0$ is
$s(s+1)$ on $\{(x,0) \mid x \neq 0\}$, and $s(s+1)(s-1)$ 
at $(0,0)$. 

Now let us compute the $b$-function for $(1-y)^su$.
The local indicial polynomial of $\delta(t+y-1)\otimes u$ 
along $t=0$ is $s(s+3)$ at any point of $t=0$.
Hence the $b$-function for $(1-y)^s u$ divides $(s+1)(s-2)$.
In the same way, the local indicial polynomial of 
$\delta(t+x-1)\otimes u$ 
along $t=0$ is $s(s+1)$ at any point of $t=0$.
Finally, the indicial polynomial of $\delta(t-x+y)\otimes u$ 
is $s(s-1)$ on $\{(x,x) \mid x \neq 0\}$, and $s(s-1)(s-2)$ 
at $(0,0)$. 

Therefore, we conclude that any rational solution to $I$, 
if it exists, can be written in the form
$g(x,y)y^{-1}(1-x)^{-1}(1-y)^{-3}$ with a polynomial $g$. 
Now we may compute the twisted ideal $I_{(0,1,0,1,3)}$,
where $f_1 = x, f_2 = y, f_3 = x-y, f_4 = x-1, f_5 = y-1$,
and $f$ is the product.  Multiplying by $f^2$, we get the expressions,
$$\begin{array}{ccl}
f^2L_1 & = & (x^2-x^3)(f\partial_x)^2 + x((1-3x)f-(1-x)y
\frac{\partial f}{\partial y} - (1-x)x\frac{\partial f}{\partial x})
(f\partial_x) + \\
& & x(1-x)y(f\partial_y)(f\partial_x) + xyf(f\partial_y) + 3xf^2 \\
f^2L_2 & = & (y^2-y^3)(f\partial_y)^2 + y((1-5y)f-(1-y)x
\frac{\partial f}{\partial x} - (1-y)y\frac{\partial f}{\partial y})
(f\partial_y) + \\
& & y(1-y)x(f\partial_x)(f\partial_y) - yxf(f\partial_x) - 3yf^2,
\end{array}$$
and we set $T_1$ and $T_2$ to be the operators obtained
from the substitution~(\ref{twistsub}).
We remark that the ideal generated by $T_1$ and $T_2$ is neither
holonomic nor specializable with respect to the weight
vector $(1,1,-1,-1)$, so it is difficult
to apply Gr\"obner deformations to it just yet.

The twisted ideal $I_{(0,1,0,1,3)}$ is the Weyl closure of the ideal 
generated by $T_1$ and $T_2$.  
At the moment, this Weyl closure
is computationally too intensive to compute.  However, we are able to
compute a partial closure by noting that $T_1$ is divisible by 
$g = f_1^3f_2f_3^2f_4f_5$ and $T_2$ is divisible by 
$h = f_1^2f_2^2f_3^2f_4f_5$.
Now the ideal $J$ generated by 
$$\begin{array}{ccl}
\frac{1}{g}T_1 & = &
(-x^3y+x^3+2x^2y-2x^2-xy+x)\partial_x^2+ \\
& &
(-x^2y^2+x^2y+2xy^2-2xy-y^2+y)\partial_x\partial_y +\\
& & (3x^2y-6xy+3y)\partial_x
+(2xy^2-2xy-2y^2+2y)\partial_y+ \\
& &
(-4xy-2x+4y+2) \\
\frac{1}{h}T_2 & = &
(-xy^4+2xy^3+y^4-xy^2-2y^3+y^2)\partial_y^2 + \\
& &
(-x^2y^3+2x^2y^2+xy^3-x^2y-2xy^2+xy)\partial_x\partial_y + \\
& & 
(3x^2y^2-4x^2y-3xy^2+x^2+4xy-x)\partial_x + \\
& &
(4xy^3-6xy^2-3y^3+2xy+4y^2-y)\partial_y + \\
& &
(-6xy^2+8xy+3y^2-2x-4y+1),
\end{array}$$
is indeed holonomic, hence we may
apply Algorithm~\ref{alg-polygrob}.
We find that the $b$-function with respect to the weight
$w = (-1,-1)$ is $(s+5)(s+2)^2(s+1)s(s-2)(s-3)^2$, which implies
that a polynomial solution to $J$ must have degree less than or equal
to $5$.  We find that $J$ has 2 polynomial solutions,
so that $I$ has 2 rational solutions,
$ (x y^2-3 x y+3 x-1)/(y-1)^3 $ and $(x-y)/y(y-1)^3$.
\end{example}

\section{Solutions by duality}

For holonomic $M$ and $N$, it is well known~\cite{KashiwaraHolonomicII} that
\begin{equation}
\label{basicisom}
\ext^i_D(M,N) \simeq 
  H^{i-n}(\Omega \otimes_D^L ({\bf D}(M) \otimes_{\k[x]}^L N)),
\end{equation}
where $\Omega := (D/\{x_1,\dots,x_n\}\cdot D)$,
and ${\bf D}(M)$ is the holonomic dual,
$$ {\bf D}(M) := \hom_{\k[x]}(\Omega, \ext^n_D(M,D)).$$

The spaces $\ext^i_D(M,N)$ are finite-dimensional ${\bf k}$-vector spaces 
and correspond to the solutions of $M$ in $N$ when $i=0$.
For example, if $N = {\bf k}[{\bf x}]$, then
we obtain the polynomial solutions of $M$, whereas if
$N = D/D\cdot\{x_1,\dots,x_n\}$, then we obtain
the delta function solutions of $M$ with support at the origin.

In this section, we explain how (\ref{basicisom}) can be
used to compute the dimensions of $\ext_D^i(M,N)$.
We first discuss how to compute the holonomic dual, next
discuss the special cases $N = {\bf C}[{\bf x}]$
and $N = {\bf C}[{\bf x}][\frac{1}{f}]$, and last discuss
the general case of holonomic $N$.
A method to extend these algorithms to compute an explicit 
basis of $\hom_D(M,N)$ and $\ext_D^i(M,N)$
is the subject of the forthcoming paper~\cite{TW}.

\bigskip
\noindent \textbf{Notation}:
Let us explain the notation we will use to write maps
of left or right $D$-modules.  As usual, maps between finitely generated
modules will be represented by matrices,
but some care has to be given to the order in which elements are
multiplied due to the noncommutativity of $D$.

Given an $r \times s$ matrix 
$A = [a_{ij}]$ with entries in $D$, we get a map
of free left $D$-modules,
$$D^r \stackrel{\cdot A}{\longrightarrow} D^s
\hspace{.2in}
[g_1,\dots, g_r] \mapsto [g_1,\dots,g_r]\cdot A,$$
where $D^r$ and $D^s$ are regarded as modules of row 
vectors, and the map is matrix multiplication.  
Under this convention, the composition of maps
$D^r \stackrel{\cdot A}{\longrightarrow} D^s$ and
$D^s \stackrel{\cdot B}{\longrightarrow} D^t$ is the map
$D^r \stackrel{\cdot AB}{\longrightarrow} D^t$ where $AB$ is 
usual matrix multiplication.
In general, suppose $M$ and $N$ are left $D$-modules 
with presentations $D^r/M_0$ and $D^s/N_0$.  Then the matrix $A$ induces
a left $D$-module map between $M$ and $N$, denoted
$(D^r/M_0) \stackrel{\cdot A}{\longrightarrow} (D^s/N_0)$,
precisely when $\vec{g}\cdot A \in N_0$ for all row vectors $\vec{g} \in M_0$.
Conversely, any map of left $D$-modules between $M$ and $N$
can be represented by some matrix $A$ in the manner above.

Now let us discuss maps of right $D$-modules.
The $r \times s$ matrix $A$ also defines a map of right $D$-modules
in the opposite direction,
$$(D^s)^T \stackrel{A \cdot}{\longrightarrow} (D^r)^T
\hspace{.2in}
[h_1,\dots, h_s]^T \mapsto A \cdot [h_1,\dots,h_s]^T,$$
where the superscript-$T$ means to regard the free modules 
$(D^s)^T$ and $(D^r)^T$ as consisting of column vectors.
This map is equivalent to the map obtained by applying $\hom_D(-,D)$
to $D^r\stackrel{\cdot A}{\longrightarrow} D^s$, thus
$(D^s)^T$ may be regarded as the dual module $\hom_D(D^s,D)$.
We will suppress the superscript-$T$ when the context is clear.
As before, the matrix $A$ induces a right
$D$-module map between right $D$-modules $N'= (D^s)^T/N'_0$ and 
$M' = (D^r)^T/M'_0$  when $A \cdot \vec{g} \in M'_0$ for all 
column vectors $\vec{g} \in N'_0$.  We denote the map by
$(D^s)^T/N'_0 \stackrel{A \cdot}{\longrightarrow} 
(D^r)^T/M'_0$.

\bigskip
\noindent \textbf{Left-right correspondence and ${\bf \Omega}$}:
As is well-known, a standard use for $\Omega$ is to establish
a correspondence between the categories of left and right
$D$-modules.
The correspondence can be expressed through
the adjoint operator $\tau$, which is the algebra involution
$$\tau : D \longrightarrow D \hspace{.3in}
x^{\alpha}\partial^{\beta} \mapsto
(-\partial)^{\beta}x^{\alpha}.$$
Namely, given a left $D$-module $M \simeq D^r/M_0$, the corresponding
right $D$-module is $\Omega \otimes_{{\bf k}[{\bf x}]} M \simeq
D^r/\tau(M_0)$.
Conversely, given a right $D$-module $N \simeq D^s/N_0$,
the corresponding left $D$-module is $\hom_{{\bf k}[{\bf x}]}(\Omega,N)
\simeq D^s/\tau(N_0)$.
Similarly, given a homomorphism of left $D$-modules 
$\phi : (D^r/M_0) {\rightarrow} (D^s/N_0)$ 
defined by left multiplication by the $r \times s$ matrix $A = [a_{ij}]$,
the corresponding homomorphism of right $D$-modules
$
\tau(\phi) : (D^r/\tau(M_0)) {\rightarrow} (D^s/\tau(N_0))
$
is defined by right multiplication by
the $s \times r$ matrix $\tau(A) := [\tau(a_{ij})]^T$.

Let us explain details of the above correspondence 
for the non-specialist.
Given a left $D$-module $M$, there is a corresponding right
$D$-module $\Omega \otimes_{\k[{\bf x}]} M$ 
where the structure is given by extending the actions,
$$(w\otimes m)f = wf \otimes m \hspace{.2in}
(w\otimes m)\xi = w\xi \otimes m - w \otimes \xi m$$
for
$f\in\k[{\bf x}]$ and $\xi\in\der(\k[{\bf x}])$.
Given a presentation $D^r/M_0$ for $M$ with
generators denoted $\{e_i\}_{i=1}^r$,
then in $\Omega \otimes_{\k[{\bf x}]} M$ we have
$$(1 \otimes e_i)x^{\alpha}\partial^{\beta} = (1 \otimes x^{\alpha}e_i)
\partial^{\beta}  =  1\otimes (-\partial)^{\beta}x^{\alpha}e_i
= 1 \otimes \tau(x^{\alpha} \partial^{\beta})e_i.$$
It follows that $\Omega \otimes_{\k[{\bf x}]} M$ is generated
by $\{1\otimes e_i\}_{i=1}^r$ and
gets the presentation
$D_n^r/\tau(M_0)$.

Conversely, given a right $D$-module $N$, there is a corresponding
left $D$-module $\hom_{\k[{\bf x}]}(\Omega, N)$ where the structure is
given by extending the action,
$$(f\varphi)(w) = \varphi(w)f \hspace{.2in}
(\xi\varphi)(w) = \varphi(w\xi) - \varphi(w)\xi$$
for
$\varphi \in \hom_{\k[{\bf x}]}(\Omega, N)$,
$w \in \Omega$, $f \in \k[{\bf x}]$,
and $\xi \in \der(\k[{\bf x}])$.
A morphism $\varphi \in \hom_{\k[{\bf x}]}(\Omega, N)$ can
be identified with its image $\varphi(1) \in N$.
Since
$$\begin{array}{cl} (x^{\alpha}\partial^{\beta}\varphi)(1) & =
(x^{\alpha}(\partial^{\beta}\varphi))(1)
= (\partial^{\beta}\varphi)(1)x^{\alpha}
= \varphi(1)(-\partial)^{\beta}x^{\alpha} \\
& =  \varphi(1)\tau(x^{\alpha}\partial^{\beta}),
\end{array}$$
the morphism $x^{\alpha}\partial^{\beta}\varphi$ gets identified
with $\varphi(1)\tau(x^{\alpha}\partial^{\beta})$.
In particular, given a presentation
$D^s/N_0$ of $N$,
then $\hom_{\k[{\bf x}]}(\Omega, N)$
is generated as a left $D$-module by the morphisms $\{\varphi_i\}_{i=1}^s$
such that $\varphi_i(1) = e_i$.  By the computation above,
a relation $\sum_i e_ig_i = 0$ in $N$ corresponds
to a relation $\sum_i \tau(g_i)\varphi_i$ in $\hom_{\k[{\bf x}]}(\Omega,N)$ 
because
$(\sum_i \tau(g_i)\varphi_i)(1) = \sum_i e_i\tau(\tau(g_i))
= \sum_i e_i g_i.$  
It follows that 
$\hom_{\k[{\bf x}]}(\Omega, N)$ is generated
by $\{\varphi_i\}_{i=1}^s$ and 
gets the presentation
\begin{equation}
\label{eq:rightleft}
{\rm Hom}_{\k[{\bf x}]}(\Omega, N) \simeq D_n^s/\tau(N_0).
\end{equation}

\subsection{Holonomic dual}
Let us discuss how ${\bf D}(M)$ can be computed.

\begin{algorithm}
\label{dual}
\rm [Computing the holonomic dual]

\noindent{\sc Input}: $D^{r_0}/D\cdot\{\vec{g}_1,\dots,\vec{g}_{r_1}\}$,
a presentation of a holonomic left $D$-module $M$.
 
\noindent{\sc Output}: The holonomic dual ${\bf D}(M)$.

\begin{enumerate}
\item Compute the first $n+1$ steps of any free resolution of $M$.  
Let the $n$-th part of the resolution be
$ D^p \stackrel{\cdot P}{\longrightarrow} D^q 
       \stackrel{\cdot Q}{\longrightarrow} D^r.$
\item Dualize and apply the adjoint operator
(recall if $P = [p_{ij}]$, then $\tau(P) = [\tau(p_{ij})]^T$) to get
$ D^p \stackrel{\cdot \tau(P)}{\longleftarrow} D^q 
       \stackrel{\cdot \tau(Q)}{\longleftarrow} D^r.$
\item Return $\ker(\cdot\tau(P))/\im(\cdot \tau(Q))$.
\end{enumerate}
\end{algorithm}

\noindent \textbf{Proof}:
Let the first $n+1$ steps of a free resolution of $M$ be denoted,
$$F^{\bullet}: D^{r_{n+1}} \stackrel{\cdot P}{\longrightarrow}
D^{r_n} \stackrel{\cdot Q}{\longrightarrow} D^{r_{n-1}}
\rightarrow \cdots \rightarrow D^{r_0} \rightarrow 0.$$
Applying $\hom_D(D,-)$ yields a complex of right $D$-modules, 
$$\hom_D(D,F^{\bullet}): (D^{r_{n+1}})^T \stackrel{P \cdot}{\longleftarrow}
(D^{r_n})^T \stackrel{Q \cdot}{\longleftarrow} (D^{r_{n-1}})^T
\leftarrow \cdots \leftarrow (D^{r_0})^T \leftarrow 0,$$
and by definition,
$$\ext^n_D(M,D) \simeq 
\frac{\ker(D^{r_{n+1}} \stackrel{P \cdot}{\longleftarrow}
D^{r_n})}{\im(D^{r_n} \stackrel{Q \cdot}{\longleftarrow} D^{r_{n-1}})}.$$

Since ${\bf D}(M) = \hom_{\k[{\bf x}]}(\Omega, 
\ext_D^n(M,D))$, it only remains to determine the effect of applying
$\hom_{\k[{\bf x}]}(\Omega, -)$. 
Using the equation~(\ref{eq:rightleft}),
if $\{\vec{L}_1,\dots,\vec{L}_k\}$ are generators of 
$K = \ker(D^{r_{n+1}} \stackrel{P \cdot}{\longleftarrow} D^{r_n})$,
and $\sum_i \vec{L}_i g_i \in I = 
\im(D^{r_n} \stackrel{Q \cdot}{\longleftarrow} D^{r_{n-1}})$
is a relation, then the corresponding relation
$\sum_i \tau(g_i)\varphi_i$ in $\hom_{\k[{\bf x}]}(\Omega, \ext_D^n(M,D))$
can be realized as the relation
$\sum_i \tau(\vec{L}_i g_i) = \tau(g_i) \tau(\vec{L}_i) \in \tau(I)$.
It follows that
$${\bf D}(M) \simeq
\frac{\ker(D^{r_{n+1}} \stackrel{\cdot \tau(P)}{\longleftarrow}
D^{r_n})}{\im(D^{r_n} \stackrel{\cdot \tau(Q)}{\longleftarrow} 
D^{r_{n-1}})},$$
which is the output of step 3.
\qed

\begin{example} 
\label{appelldual}
\rm
The Appell differential equation $F_1(2,-3,-2,5)$ of Example \ref{appell}
has the resolution
$\displaystyle{0 {\rightarrow} D^1
\stackrel{\cdot Q_1}{\rightarrow} D^2 
\stackrel{\cdot Q_0}{\rightarrow} D^1
\rightarrow 0,}$
where
$$ \begin{array}{l}
Q_0 = \left[ \begin{array}{l}
(\theta_x-3)\partial_y - (\theta_y-2)\partial_x \\
(y^2-y)(\partial_x\partial_y + \partial_y^2)-2(y+x)\partial_x
+4y\partial_y+2\partial_x-8\partial_y-4 
\end{array} \right]^T \\
 \\
Q_1 =  \left[ \begin{array}{c}
(y^2-y)(\partial_x\partial_y+\partial_y^2)
-2x\partial_x+6y\partial_y+\partial_x-9\partial_y \\
-(\theta_x-3)\partial_y + (\theta_y-1)\partial_x
\end{array} \right]
\end{array}$$
The holonomic dual ${\bf D}(F_1(2,-3,-2,5))$ is the cokernel
of $\tau(Q_1)$ and is the Appell differential equation 
$F_1(-1,4,2,-3)$.

\end{example}

\subsection{Polynomial and rational solutions by duality}

When $N = \k[{\bf x}]$, the isomorphism (\ref{basicisom})
specializes to
\begin{equation}
\label{polyisom}
\ext^i_D(M,\k[{\bf x}]) \simeq 
  H^{n-i}(\Omega \otimes_D^L {\bf D}(M)).
\end{equation}
The right hand side is equivalently the $(n-i)$-th integration of
${\bf D}(M)$ to the origin.  An algorithm to compute
integration is given in~\cite{OT}.  Using it, we can evaluate the
dimensions of $\ext^i_D(M,\k[{\bf x}])$ and in particular
$\hom_D(M,\k[{\bf x}])$. 

\begin{algorithm}\rm [Evaluating dimensions of polynomial solution spaces]

\noindent{\sc Input}: a holonomic left $D$-module $M$.
 
\noindent{\sc Output}: dimensions of $\ext_D^i(M,\k[{\bf x}])$.

\begin{enumerate}
\item Compute the dual ${\bf D}(M)$ using Algorithm \ref{dual}
\item Compute the integrations of ${\bf D}(M)$ to the origin using the 
algorithm in~\cite{OT}.  They are finite dimensional
vector spaces.
\item Return the dimensions.
\end{enumerate}
\end{algorithm}

The dimensions of rational solution spaces can be evaluated in a similar
way.  When $N = \k[{\bf x}][\frac{1}{f}]$,
the isomorphism (\ref{basicisom}) specializes to
\begin{equation}
\label{ratlisom}
\ext^i_D(M,\k[{\bf x}][{1}/{f}]) \simeq 
  H^{n-i}(\Omega \otimes_D^L {\bf D}(M)[{1}/{f}]).
\end{equation}
The right hand side is now equivalently the $(n-i)$-th integration of
${\bf D}(M)[\frac{1}{f}]$ to the origin.  An algorithm to compute
localization is given in~\cite{OTW}.  Using it and the integration algorithm, 
we can evaluate the dimensions of $\ext^i_D(M,\k[{\bf x}][\frac{1}{f}])$ and
$\hom_D(M,\k[{\bf x}][\frac{1}{f}])$. To get the dimension
of all rational solutions, take $f$ to be any polynomial vanishing on
the singular locus.

We summarize how to compute 
the integration of a module $N$ to the origin according
to~\cite{OT} in a slightly more general way.  
The generalization sometimes gives a more efficient strategy
than~\cite{OT}. We need
to recall some definitions.  To any
strictly positive $w \in {\bf Z}^n_{>0}$, 
we get an integration filtration $F_w$ of $D$ defined by
$F_w^i(D) = \span_k \{x^{\alpha}\partial^{\beta} |
w\cdot \alpha - w\cdot \beta \leq i\}$.
More generally, for $\vec{m} \in Z^{r}$,
we also get a shifted filtration $F_w[\vec{m}]$ of the free module
$D^r$ defined by $F_w^i[\vec{m}](D^r) = 
\span_k \{x^{\alpha}\partial^{\beta} e_j | 
w\cdot \alpha - w\cdot \beta - m_j\leq i\}$.
We will often write $D^r[\vec{m}]$ for the free module $D^r$ equipped
with the shifted filtration $F_w[\vec{m}]$ when the context is clear.
The filtrations $F_w[\vec{m}]$ induce filtrations on subquotients
of $D^r$ in the natural way.
Now we may say the steps of the integration algorithm.
First, compute a $(w,-w)$-strict free resolution $G^{\bullet}$ of $N$
of length $n+1$.  This is a resolution of
$N$ by free modules $D^{r_j}[\vec{m}_j]$ with the property that
the differentials preserve the filtration and moreover
induce a resolution on the associated graded level.
Second, compute
the integration $b$-function of $N$ with respect to $(w,-w)$, and find
its minimal and maximal integral roots $k_0$ and $k_1$.
The integration $b$-function is the monic polynomial
$b(s)$ of least degree satisfying $b( \sum_i w_i\partial_i x_i )
\cdot F^0(N) \subset F^{-1}(N)$.
Third, compute the cohomology of the complex
${F_{-k_0}(\Omega \otimes_D G^{\bullet})}/
{F_{-k_1-1}(\Omega \otimes_D G^{\bullet})}$,
which is a complex of finite-dimensional vector spaces.  The dimensions
of the cohomology groups are equal to the dimensions of the integration
modules of $N$.

\begin{example} \rm
Let us evaluate the dimension of polynomial solutions
to the Appell differential equation $M = F_1(2,-3,-2,5)$ of Example 
\ref{appell}. Choose the weight 
vector $w = (1,2)$.  The resolution of Example~\ref{appelldual}, after
dualizing, applying the adjoint operator, and shifting,
$$0 {\rightarrow} D^1[0]
\stackrel{\cdot \tau(Q_0)}{\longrightarrow} D^2[-1,1] 
\stackrel{\cdot \tau(Q_1)}{\longrightarrow} D^1[0]
\rightarrow 0,$$
preserves filtrations but does not induce a resolution
on the associated graded level.  On the other hand,
if we adjust the resolution to
$$G^{\bullet}: 0 {\rightarrow} D^1[1]
\stackrel{\cdot P_0}{\longrightarrow} D^2[0,1] 
\stackrel{\cdot P_1}{\longrightarrow} D^1[0]
\rightarrow 0,$$
where
$$P_0 = \left[
\begin{array}{c}
-(\theta_x + 5)\partial_y + (\theta_y + 2)\partial_x \\
(x^2-x)(\partial_x^2+\partial_x \partial_y) + 
4x\partial_x + 2(3x+2y)\partial_y + 4\partial_x - 5\partial_y - 2
\end{array}
\right]^T$$
$$P_1 = \left[
\begin{array}{c}
(x^2-x)(\partial_x^2+\partial_x \partial_y) + 
2x\partial_x + 4(x+y)\partial_y +5\partial_x - 4\partial_y - 4\\
(\theta_x + 4)\partial_y - (\theta_y + 2)\partial_x
\end{array}
\right],$$
then we do obtain a $(w,-w)$-strict resolution of
${\bf D}(M) = F_1(-1,4,2,-3)$.

The integration $b$-function with respect to $(w,-w)$ 
is $(s+5)(s-2)(s-5)$, hence the integration complex for ${\bf D}(M)$
is quasi-isomorphic to the truncated complex
$F_w^5(\Omega \otimes_D G^{\bullet}) / F_w^{-6}(\Omega \otimes_D G^{\bullet})$,
which is a complex of finite-dimensional vector spaces with dimensions,
$$0 \rightarrow {\bf Q}^{16} \stackrel{\cdot P_0}{\longrightarrow} 
{\bf Q}^{28} \stackrel{\cdot P_1}{\longrightarrow} 
{\bf Q}^{12} \rightarrow 0.$$
For instance, $F_w^5(\Omega[0])$ consists of the $12$ monomials,
$$\{1,x,y,x^2,xy,y^2,x^3,x^2y,xy^2,x^4, xy^3, x^5\},$$ 
and so on.  Note that $\tau(P_1)$ is a $(w,-w)$-Gr\"obner basis
of $F_1(2,-3,-2,5)$ and hence for this case, the
duality method essentially coincides with the Gr\"obner 
deformation method of Section \ref{GBdeform} at the level
of ${\rm Hom}_D(M, \k[{\bf x}])$.
The above computations were made in Macaulay 2, where we get
the output,
$$
\begin{array}{cclllll}
{\tt \scriptstyle i1} & {\tt \scriptstyle :} & 
{\tt \scriptstyle PolyExt(M)} & & \\
{\tt \scriptstyle o1} & {\tt \scriptstyle =} & {\tt \scriptstyle HashTable\{} & {\tt \scriptstyle 0} & {\tt \scriptstyle =>} &
{\tt \scriptstyle QQ^1} & \}  \\
& & & {\tt \scriptstyle 1} & {\tt \scriptstyle =>} & {\tt \scriptstyle QQ^2} &   \\
& & & {\tt \scriptstyle 2} & {\tt \scriptstyle =>} & {\tt \scriptstyle QQ^1} &   
\end{array}
$$
Here, the output ${\tt \scriptstyle i\ =>\ QQ^j}$ means that  
$\dim \ext_D^i(M,\k[{\bf x}]) = j$.

\end{example}

\begin{example} \rm  Let us now evaluate the dimension
of rational solutions to $M = F_1(2,-3,-2,5)$.
The singular locus is $xy(x-y)(x-1)(y-1)$.  
We will search for solutions 
in $\k[x,y][\frac{1}{x}]$ first.  
From Example \ref{appelldual},
${\bf D}(M)$ has the presentation $D/\tau(Q_1)$.  Let $u$ be the section
corresponding to the residue class of $\overline{1}$ in this presentation.
Then the localization ${\bf D}(M)[\frac{1}{x}]$
is generated by $u \otimes \frac{1}{x^7}$ and gets the presentation
$D/J$, where
$$J = D \cdot \left\{ \begin{array}{c}
(\theta_x\theta_y + \theta_y^2 + 8\theta_y + 2\theta_x + 12)
-(\theta_x+\theta_y+4)\partial_y \\
(\theta_x\theta_y+2\theta_x+7\theta_y+14) - (\theta_x+10)x\partial_y
\end{array} \right\}.$$
The natural localization map can be written
as $\varphi : D/\tau(Q_1) \longrightarrow D/J$, where $\varphi(1) = x^7$.
Choose the integration weight vector $w = (1,2)$.  
Then ${\bf D}(M)[\frac{1}{x}]$ has a $(w,-w)$-strict resolution
$$G^{\bullet}: 0 \rightarrow D^1[-1] \stackrel{\cdot[v_1,v_2]}
{\longrightarrow} D^2[0,-1]
\stackrel{\cdot[u_1,u_2]^T}{\longrightarrow} D^1[0] \rightarrow 0$$
where
$$\begin{array}{ccl}
u_1 & = & -x^2\partial_x\partial_y+xy\partial_x\partial_y+2x\partial_x-11x\partial_y+7y\partial_y+14 \\
u_2 & = & x^3\partial_x^2+x^3\partial_x\partial_y-x^2\partial_x^2-x^2\partial_x\partial_y+16x^2\partial_x+11x^2\partial_y + \\
& & 4xy\partial_y-9x\partial_x-11x\partial_y+52x-7 \\
v_1 & = & x^3\partial_x^2+x^3\partial_x\partial_y-x^2\partial_x^2-x^2\partial_x\partial_y+16x^2\partial_x+12x^2\partial_y+ \\
& & 4xy\partial_y-8x\partial_x-11x\partial_y+52x-6 \\
v_2 & = & x^2\partial_x\partial_y-xy\partial_x\partial_y-2x\partial_x+11x\partial_y-6y\partial_y-12
\end{array}$$
The integration $b$-function is $(s+12)(s+5)(s+2)$, hence we want
the cohomology of the complex
$F_w^{12}(\Omega \otimes_D G^{\bullet})/F_w^1(\Omega \otimes_D G^{\bullet})$,
which has the shape,
$$0 \rightarrow {\bf Q}^{41} \longrightarrow {\bf Q}^{86}
\longrightarrow {\bf Q}^{45} \rightarrow 0.$$
By evaluating the dimensions of the cohomology groups
in Macaulay 2, we find
$$
\begin{array}{cclllll}
{\tt \scriptstyle i2} & {\tt \scriptstyle :} & {\tt \scriptstyle RatlExt(M, x)} & & \\
{\tt \scriptstyle o2} & {\tt \scriptstyle =} & {\tt \scriptstyle HashTable\{} & {\tt \scriptstyle 0} & {\tt \scriptstyle =>} &
{\tt \scriptstyle QQ^2} & \}  \\
& & & {\tt \scriptstyle 1} & {\tt \scriptstyle =>} & {\tt \scriptstyle QQ^5} &   \\
& & & {\tt \scriptstyle 2} & {\tt \scriptstyle =>} & {\tt \scriptstyle QQ^3} &   
\end{array}
$$
Since we already computed a polynomial solution, this means there is one
rational solution with pole along $x$.  Similarly, 
we get the exact same dimensions for
$\ext_D^i(M,\k[x,y][{\scriptstyle \frac{1}{y}}])$,
which means that there is also one rational solution with pole along $y$.
The rank of the system is $3$, therefore we have found all the solutions.
We could also compute,
$$
\begin{array}{cclllll}
{\tt \scriptstyle i3} & {\tt \scriptstyle :} & {\tt \scriptstyle RatlExt(M, f)} & & \\
{\tt \scriptstyle o3} & {\tt \scriptstyle =} & {\tt \scriptstyle HashTable\{} & {\tt \scriptstyle 0} & {\tt \scriptstyle =>} &
{\tt \scriptstyle QQ^1} & \}  \\
& & & {\tt \scriptstyle 1} & {\tt \scriptstyle =>} & {\tt \scriptstyle QQ^3} &   \\
& & & {\tt \scriptstyle 2} & {\tt \scriptstyle =>} & {\tt \scriptstyle QQ^2} &   
\end{array}
$$
where $f$ is any of the polynomials $x-y$, $x-1$, or $y-1$.
As expected, there are no rational solutions with poles along $x-y$,
$x-1$, or $y-1$, but in all cases 
there are new $\ext^1$ and $\ext^2$.
We have not computed $\ext$ with respect to any products of poles since
it is computationally too intensive for now.
\end{example}

Once we have evaluated the dimension of the solution spaces,
we can compute the solutions by a brute force method.
\begin{enumerate}
\item For a given holonomic system $M$, compute its singular locus.
Let $f$ be a polynomial such that $f=0$ contains the singular locus.
\item Evaluate the dimension $d$ of the rational solutions 
by the homological duality method.
\item Try to find rational solutions of the form $r(x)/f^k$, 
${\rm degree}(r) = p$.
Increase $p+k$ until we find $d$ linearly independent solutions.
\end{enumerate}

\subsection{Holonomic solutions by duality}
The isomorphism (\ref{basicisom}) can also be expressed 
(see e.g.~\cite{Bj}) as
$$\ext_D^i(M,N) \simeq \tor_{n-i}^D(\ext^n_D(M,D), N).$$
To compute the right-hand-side, it is also well known that for 
all $M'$, and in particular $M'=\ext^n_D(M,D)$, 
$$\tor^D_{n-i}(M',N) \simeq H^{i}(K^{\bullet}( (M'
\otimes_{\k[\bf x]} \Omega^{-1}) 
\widehat{\otimes} N ;
\{x_i - y_i, \partial_i +\delta_i\}_{i=1}^n)),$$
where  
$K^{\bullet}$ denotes the Koszul complex
and $\widehat{\otimes}$ denotes the external tensor product
into the category of $D_{2n} = {\bf k}\langle
x_1,y_1,\dots,x_n,y_n,\partial_1,\delta_1,\dots,\partial_n,
\delta_n\rangle$-modules.  Combining these
isomorphisms leads to
$$\ext_D^i(M,N) \simeq
H^{i}(K^{\bullet}( ({\bf D}(M)
\widehat{\otimes} N ;
\{x_i - y_i, \partial_i +\delta_i\}_{i=1}^n))$$
By an automorphism of $D$, we can transform
$\{x_i - y_i, \partial_i +\delta_i\}_{i=1}^n$ into
$\{x_i,y_i \}_{i=1}^n$, for which the Koszul complex
computes the derived restriction to the origin.

\begin{algorithm}\rm [Evaluating dimensions of holonomic solution spaces]

\noindent{\sc Input}: holonomic left $D$-modules $M$ and $N$
 
\noindent{\sc Output}: dimensions of $\ext_D^i(M,N)$.

\begin{enumerate}
\item Compute the dual ${\bf D}(M)$ using Algorithm \ref{dual}
\item Form the $D_{2n}$-module ${\bf D}(M) \otimes_k N$
and apply the change of coordinates
$\eta: D_{2n} \rightarrow D_{2n}$ where $\eta$ maps,
$$\begin{array}{lcl}
x_i \mapsto \frac{1}{2}x_i- \delta_i, & & 
\partial_i \mapsto \frac{1}{2}y_i+\partial_i, \\
y_i \mapsto -\frac{1}{2}x_i- \delta_i, 
& & \delta_i 
\mapsto \frac{1}{2}y_i- \partial_i \end{array}.$$
\item Compute the restrictions of $\eta({\bf D}(M) \otimes_k N)$
to the origin using the algorithm in~\cite{OT}.  They are finite dimensional
vector spaces.
\item Return the dimensions.
\end{enumerate}
\end{algorithm}

\begin{example} \rm
Let $M = F_1(2,-3,-2,5)$ be the Appell differential equation
of example \ref{appell}, and let
$N = \k[x,y][\frac{1}{x}]/\k[x,y]$. It has presentation $D/D\cdot
\{x,\partial_y\}$, where the generator $1$ corresponds to $\frac{1}{x}$.
Using the above algorithm, we compute
$$
\begin{array}{cclllll}
{\tt \scriptstyle i4} & {\tt \scriptstyle :} & {\tt \scriptstyle DExt(M,N)} & & \\
{\tt \scriptstyle o4} & {\tt \scriptstyle =} & {\tt \scriptstyle HashTable\{} & {\tt \scriptstyle 0} & {\tt \scriptstyle =>} &
{\tt \scriptstyle QQ^1} & \}  \\
& & & {\tt \scriptstyle 1} & {\tt \scriptstyle =>} & {\tt \scriptstyle QQ^3} &   \\
& & & {\tt \scriptstyle 2} & {\tt \scriptstyle =>} & {\tt \scriptstyle QQ^2} &   
\end{array}
$$

Similarly, let $N = \k[\partial_x, \partial_y] \simeq
D/D\cdot\{x,y\}$, the module of the delta functions with the support $(0,0)$.  
Then we compute
$$
\begin{array}{cclllll}
{\tt \scriptstyle i5} & {\tt \scriptstyle :} & {\tt \scriptstyle DExt(M,N)} & & \\
{\tt \scriptstyle o5} & {\tt \scriptstyle =} & {\tt \scriptstyle HashTable\{} & {\tt \scriptstyle 0} & {\tt \scriptstyle =>} &
{\tt \scriptstyle QQ^0} & \}  \\
& & & {\tt \scriptstyle 1} & {\tt \scriptstyle =>} & {\tt \scriptstyle QQ^1} &   \\
& & & {\tt \scriptstyle 2} & {\tt \scriptstyle =>} & {\tt \scriptstyle QQ^2} &   
\end{array}
$$
 \end{example}

As before, once we know the dimension of $\hom_D(M,N)$, we
can compute the solutions of $M$ in $N$ by a brute force
method.

\begin{enumerate}
\item For given holonomic systems $M$ and $N$, 
evaluate the dimension $d$ of $\hom_D(M,N)$
by the homological duality method.
\item Filter $N$ by finite-dimensional vector spaces $F^i(N)$
and search for solutions in $F_i(N)$ for increasing $i$ until
$d$ linearly independent solutions are found.  
\end{enumerate}

For instance in step 2, if $N = D/J$, then we can use the induced
Bernstein filtration $B$ where $B^i(D/J)$ consists of
residues of elements $L \in D$ whose total degree is less than
or equal to $i$.

\rightline{ Toshinori Oaku, \quad oaku@twcu.ac.jp }
\rightline{ Department of Mathematics, Tokyo Woman's Christian University}
\rightline{ 2-6-1 Zempukuji, Suginami-ku, Tokyo 167-8585, Japan }

\medskip

\rightline{ Nobuki Takayama, \quad takayama@math.kobe-u.ac.jp}
\rightline{ Department of Mathematics, Kobe University}
\rightline{ Kobe  657-8501, Japan}

\medskip

\rightline{ Harrison Tsai, \quad htsai@math.berkeley.edu}
\rightline{ Department of Mathematics, University of California, Berkeley}
\rightline{ Berkeley, CA 94720-3840, USA}

\end{document}